\documentclass[11pt]{article}
\usepackage{fullpage,url,natbib}
\usepackage{amsthm,amsfonts,amsmath,amssymb,epsfig,color,float,graphicx,verbatim}
\usepackage{hyperref}
\usepackage{tikz}
\hypersetup{
	colorlinks   = true, %
	urlcolor     = blue, %
	linkcolor    = blue, %
	citecolor   = blue %
}

\newtheorem{theorem}{Theorem}
\newtheorem{proposition}{Proposition}
\newtheorem*{proposition*}{Proposition}
\newtheorem{lemma}{Lemma}

\renewcommand{\eqref}[1]{Eq.~(\ref{#1})}

\newcommand{\R}{\mathbb{R}}

\newcommand{\N}{\mathbb{N}}

\newcommand{\inv}{^{-1}}

\DeclareMathOperator*{\argmin}{arg\,min}

\newcommand{\abs}[1]{\left| #1 \right|}
\newcommand{\paren}[1]{\left( #1 \right)}
\newcommand{\sqprn}[1]{\left[ #1 \right]}
\newcommand{\nrm}[1]{\left\Vert #1 \right\Vert}
\newcommand{\norm}[1]{\Vert #1 \Vert}

\newcommand{\vertiii}[1]{{\left\vert\kern-0.25ex\left\vert\kern-0.25ex\left\vert #1 
    \right\vert\kern-0.25ex\right\vert\kern-0.25ex\right\vert}}

\newcommand{\mexp}{\mathbb{E}}

\newcommand{\PR}[2][]{\mathop{\mathbb{P}}_{#1}\left( #2 \right)}
\newcommand{\E}{\mathop{\mexp}}
\renewcommand{\P}{\mathbb{P}}

\newcommand{\eps}{\varepsilon}

\newcommand{\eqdef}{:=}

\usepackage{xcolor}

\newcommand{\supr}[1]{^{(#1)}}

\def\longto{\mathop{\longrightarrow}\limits}

\newcommand{\ninf}{\longto_{n\to\infty}}

\newcommand{\mx}{\vee}
\newcommand{\mn}{\wedge}

\newcommand{\calP}{\mathcal{P}}
\newcommand{\calQ}{\mathcal{Q}}

\newcommand{\lgc}{\mathsf{LGC}}

\newcommand{\p}[1][]{p_{#1}}
\newcommand{\minp}[1][]{\dot{p}_{#1}}
\newcommand{\pn}[1][]{\hat p_{#1}}
\newcommand{\tpn}[1][]{\tilde p_{#1}}

\newcommand{\pp}{
[0,
{\textstyle\frac12}
]^\N_{\downarrow0}
}

\newcommand{\ppoh}{
[0,
{\textstyle\frac12}
]^\N
}

\newcommand{\decr}[1]{{#1}^{\downarrow}}

\newcommand{\beq}{\begin{eqnarray*}}
\newcommand{\eeq}{\end{eqnarray*}}
\newcommand{\beqn}{\begin{eqnarray}}
\newcommand{\eeqn}{\end{eqnarray}}

\newcommand{\QED}{\hfill\ensuremath{\square}}

\newcommand{\ceil}[1]{{\left\lceil#1\right\rceil}}
\newcommand{\floor}[1]{{\left\lfloor#1\right\rfloor}}

\usepackage{ifmtarg}
\usepackage{xifthen}%
\newcommand{\ent}[1][]{%
\ifthenelse{\isempty{#1}}{%
\mathrm{H}
}{
\mathrm{H}^{(#1)}
}}

\newcommand{\loch}[1][]{%
\ifthenelse{\isempty{#1}}{%
\mathrm{h}
}{
\mathrm{h}^{(#1)}
}}

\newcommand{\mathe}{\mathrm{e}}

\newcommand{\hide}[1]{}

\newcommand{\set}[1]{\left\{ #1 \right\}}

\newcommand{\Bernu}{\operatorname{Bernoulli}}

\newcommand{\Binom}{\operatorname{Binomial}}

\newtheorem*{rep@theorem}{\rep@title}
\newcommand{\newreptheorem}[2]{%
\newenvironment{rep#1}[1]{%
 \def\rep@title{#2 \ref{##1}}%
 \begin{rep@theorem}}%
 {\end{rep@theorem}}}
\makeatother

\newreptheorem{proposition}{Proposition}

\renewcommand{\eqref}[1]{(\ref{#1})}

\title{
The Empirical Mean is Minimax Optimal for Local Glivenko-Cantelli
}

\author{%
  Doron Cohen \\
  \texttt{doronv@post.bgu.ac.il}
\and
  Aryeh Kontorovich \\
  \texttt{karyeh@cs.bgu.ac.il}
\and
Roi Weiss \\
  \texttt{roiw@ariel.ac.il }
  }

\date{}

\begin{document}

\maketitle
\begin{center}\vspace{-1cm}\today\vspace{0.5cm}\end{center}

\begin{abstract}
We revisit the recently introduced Local Glivenko-Cantelli setting, which studies distribution-dependent uniform convergence rates of the Empirical Mean Estimator (EME). In this work, we investigate generalizations of this setting where arbitrary estimators are allowed rather than just the EME. Can a strictly larger class of measures be learned? Can better risk decay rates be obtained? We provide exhaustive answers to these questions—which are both negative, provided the learner is barred from exploiting some infinite-dimensional pathologies. On the other hand, allowing such exploits does lead to a strictly larger class of learnable measures.
\end{abstract}

\section{Introduction}

\citet{CohenK23} 
initiated the study of 
the {\em local Glivenko-Cantelli}
setting: laws of large numbers that are uniform
over a function class but rather than being
{\em universal} over all distributions,
feature a delicate dependence of the risk decay on
the (local) sampling measure.
This naturally led to the 
{\em binomial empirical process}: 
for a fixed $p\in[0,1]^\N$
and
each $n\in\N$,
we have a sequence
of independent 
$Y_j\sim\Binom(n,\p[j])$,
which are centered and normalized
to obtain
$\bar Y_j:=
n\inv Y_j-\p[j]
$.
The object of interest is the
expected uniform absolute deviation:
\beqn
\label{eq:barY}
\Delta_n &:=&
\E\sup_{j\in\N}|\bar Y_j|.
\eeqn

More generally, one could imagine
fixing a distribution $\mu$ on $\set{0,1}^\N$,
sampling $X\supr1,X\supr2,\ldots,X\supr n$
i.i.d.\ from $\mu$, and estimating $p:=\E X\supr1$
via the Maximum-Likelihood Estimator (MLE)
$\pn:=\frac1n\sum_{i=1}^n X\supr i$.
In the case where $\mu$ is a product measure
(that is, the components of the vector $X\sim\mu$ are mutually independent), 
$
\E\nrm{\pn-p}_\infty$
recovers
the expression in \eqref{eq:barY}.
Despite its austere appearance,
the binomial empirical process with independent
coordinates $Y_j$ under $\ell_\infty$-norm
deviation already captures much of the richness of problem. Extensions to more general 
product
distributions
$\mu$
over $[0,1]^\N$
are straightforward
\citep[Corollary 6]{blanchard2023tight}
and the behavior under $\ell_r$
norms for $r<\infty$
is considerably simpler (Proposition 7 ibid.).
Finally, the in-expectation bounds are readily
converted to high-probability tail bounds
(Proposition 9 ibid.), and all of the upper
bounds stated for product measures hold verbatim
for arbitrary correlations.

For the purpose of analyzing \eqref{eq:barY},
\citeauthor{CohenK23}
showed that there is no loss of generality in
restricting $p$ to the set $\pp$,
consisting of all $p\in\ppoh$
with $\p[j]\downarrow 0$.
They defined $\lgc\subset\pp$
as the family of $p$
for which 
$
\Delta_n\ninf0
$
and showed that
$\lgc$
consists of exactly those $p
$
for which
\beqn
\label{eq:Told}
T(p) &:=& \sup_{j\in\N}\frac{\log (j+1)}{\log(1/\p[j])},
\qquad
p\in\pp
\eeqn
is finite.
They also characterized up to constants the asymptotic decay of $\Delta_n$ 
(whenever $T(p)<\infty$) via the functional
\beqn
\label{eq:Sold}
S(p) &:=& \sup_{j\in\N}\p[j]\log (j+1),
\qquad
p\in\pp,
\eeqn
establishing that $\Delta_n(p)$ decays as $\sqrt{S(p)/n}$.  Additional finite-sample bounds provided therein were tightened by \citet{blanchard2023tight} as follows:
\begin{align*}
    \Delta_n(p) \asymp & 1\land \paren{\sqrt{\frac{S(p)}{n}} + \sup_{j\geq 1}\frac{\log(j+1)}{n\log\paren{2+\frac{\log(j+1)}{n\p[j]}}}}, 
    \\&\text{if } n\cdot \sup_{j\geq 1} 2j\p[j] >1, 
    \\
    \Delta_n(p) \asymp & \frac{1}{n}\land \sum_{j\geq 1}\p[j], 
    \\& \text{otherwise.} 
\end{align*}
In a later work, 
\citet{BlanchardCK24}
extended some of the analysis to the much more
difficult case where $\mu$
is not a product measure (i.e.,
the coordinates of
$X\sim\mu$
have correlations).
In the present paper,
we return to the product-measure case
and investigate a different extension:
How does $\Delta_n$ behave if rather than
restricting the estimator to the MLE $\pn$,
we allow {\em arbitrary} estimators $\tilde p$?

Actually, a bit of a refinement in terminology is necessary. When considering 
(essentially) unrestricted classes of distributions such as those parametrized by
$p\in\pp$, the Empirical Mean Estimator (EME)
and
the Maximum Likelihood Estimator (MLE)
coincide. However, for more general
families $\calP\subset[0,1]^\N$, this is no longer the case: the likelihood of a given sample might be maximized over $\calP$
by some $\hat p$ other than the EME.
Hence, in the sequel, we shall be pedantic about this distinction, focusing on the EME as the more natural candidate.

Formally, an {\em estimator}
$\tilde p$
is 
any mapping from $(\set{0,1}^\N)^n$ to $[0,1]^\N$.
Any $p\in[0,1]^\N$ induces the product measure
\beqn
\label{eq:mup}
\mu = \mu(p) = \Bernu(\p[1])\otimes \Bernu(\p[2]) 
\otimes \ldots
\eeqn
on $\set{0,1}^\N$.
If 
$X\supr1,X\supr2,\ldots,X\supr n$
are sampled
i.i.d.\ from $\mu$,
then these induce
$\tilde\Delta_n:=
\E
\nrm{\tilde p-p}_\infty
$.
We say that a family of product
distributions induced by
$\calP\subset[0,1]^\N$
is {\em learnable by $\tilde p$}
if $\tilde\Delta_n\ninf0$
for each $p\in\calP$,
and just {\em learnable}
if it is learnable by some $\tilde p$.
(Since the sequence $p$
fully determines the measure $\mu(p)$,
it is fitting to say that $\tilde p$
``learns'' $p$ --- and hence also $\mu(p)$.)

This general setting
immediately raises the natural
questions: Can 
$\lgc$ be expanded to a larger learnable family via some
estimator $\tilde p$ different from the EME? Can 
some estimator $\tilde p$ achieve better decay rates
for $\tilde \Delta_n$ than the EME?

\paragraph{Our contributions.}
Modulo some technical caveats, we resolve both
questions above in the negative.
If the learner is barred from exploiting
some pathological quirks of the infinite-dimensional setting, then essentially $\lgc$
as defined above is the largest learnable family
(Theorem~\ref{thm:EMELGC}).
Furthermore, the EME achieves the minimax
risk decay rate over non-pathological distribution families
(Theorem~\ref{thm:minimax}).
Finally, in
Theorem~\ref{thm:relax}
we show that 
non-trivial extensions of $\lgc$ become possible
once the restrictions are relaxed.

\paragraph{Related work.}
Estimating the mean of a 
high-dimensional
distribution 
from independent draws is among the most basic problems of statistics.
Much of the earlier theory has focused on obtaining efficient estimators $\hat m_n$ of the true mean $m$ and analyzing the decay of $\nrm{\hat m_n-m}_2$ as a function of sample size $n$, dimension $d$, and various moment assumptions on $X$
\citep{catoni2012challenging,Devroye16,LuMen19a,LuMen19b,Cherapanamjeri19,Cherapanamjeri20,diakonikolas2020outlier,hopkins2020mean,LuMen21,lee2022optimal}. 
For $d$-dimensional distributions $\mu$ on $\{0,1\}^d$, 
Chernoff and union bounds yield
$\Delta_n(\mu )\lesssim \sqrt{\ln(d+1)/n}$
for the EME, and a simple information-theoretic argument
shows that this is minimax-optimal up to constants
\citep[Proposition 1]{CohenK23}. \citeauthor{CohenK23} further motivated
their choice of the $\ell_\infty$ norm as the most interesting of all the $\ell_r$ norms,
in a well-defined sense
(see \citet[Proposition 7]{blanchard2023tight}).
\citet{blanchard2023tight} fully closed the gaps in the analysis of \citeauthor{CohenK23}, and \citet{BlanchardCK24} took the first
nontrivial steps in analyzing non-product sampling distributions.

\paragraph{Notation.}
The measure-theoretic subtleties of defining distributions on $\set{0,1}^\N$ are addressed in \citet{CohenK23}. Our logarithms will always be base $\mathe$ by default; other bases will be explicitly specified.
The natural numbers
are denoted by
$\N=\set{1, 2, 3, \dots}$
and
for $k\in\N$, we write $[k]=\set{
i\in\N:i\le k
}$.
\hide{
For $n,d\in\N$ and a distribution $\mu$
over $\set{0,1}^d$,
we will always denote by
$p,\pn\in[0,1]^d$ the true and
empirical means of $\mu$, respectively,
as defined immediately preceding \eqref{eq:kdef}.
The $j$th coordinate of a vector $v\in\R^d$
is denoted by $v(j)$.
The expected maximal deviation
$\Delta_n(\mu)$ is defined in \eqref{eq:kdef},
and
the notation
$\pp$, $\lgc$
from the preceding paragraph
will be used throughout.
We say that $\mu$ is a {\em product}
distribution on $\set{0,1}^d$
if it can be expressed as a 
tensor
product
of $d$ distributions:
$\mu=\mu_1\otimes\mu_2\otimes\ldots\otimes\mu_d$;
equivalently, for $X\sim\mu$ we have that
the random variables $\set{X(j):j\in[d]}$
are mutually independent.}
The floor and ceiling functions, $\floor{t}$, $\ceil{t}$, map $t\in\R$ to its closest integers below and above, respectively; also, $s\mx t \eqdef \max\set{s,t}$, $s\mn t \eqdef \min\set{s,t}$.
Unspecified constants such as $c,c'$ may change value from line to line.
We use superscripts to denote distinct random vectors and subscripts to denote indices within a given vector. Thus, if $X^{(1)},,\ldots,X^{(n)}$ are independent copies of $X$, then $X^{(i)}_j$ denotes the $j$th entry of the $i$th copy.

When considering the EME as the sole estimator (as in previous works),
no generality was lost in 
restricting the range
of $p$
to $[0,\frac12]$
and
assuming sequences
monotonically decreasing to $0$
(i.e., $\pp$). The definitions of $T$ and $S$ in (\ref{eq:Told}, \ref{eq:Sold})
were based on this assumption. In this work, we will need their slightly generalized versions.
With the convention $\dot{x}:=\min\set{x,1-x}$, we define
\beqn
\label{eq:Tnew}
{T}(p) &:=& \inf_{\sigma:\N\to\N}\sup_{j\in\N}\frac{\log (j+1)}{\log(1/\minp[\sigma(j)])}, \\
\label{eq:Snew}
{S}(p) &:=& \inf_{\sigma:\N\to\N}\sup_{j\in\N}\minp[\sigma(j)]\log (j+1),
\eeqn
for $p\in[0,1]^\N$, where the infimum is over all permutations $\sigma$ over $\N$.
Whenever $\minp[j]\to0$, a unique non-increasing permutation $\decr{\minp}$
exists, and it is easily seen to be the one achieving both infima above;
thus, for $p
\in\pp
$,
the definitions in 
(\ref{eq:Tnew}, \ref{eq:Snew})
coincide with those in
(\ref{eq:Told}, \ref{eq:Sold}).

Any $p\in[0,1]^\N$ defines the product measure $\mu=\mu(p)$ as in \eqref{eq:mup}. An {\em estimator} 
$\tilde p$
and its induced deviation $\tilde \Delta_n$
are defined just above \eqref{eq:mup}, and the {\em learnability} of a
family
$\calP\subset[0,1]^n$ is defined just below
it.

We say that a family
$\calP\subset[0,1]^\N$
is {\em decaying} if 
${\displaystyle \lim_{j\to\infty}\minp[j]}=0$
for all $p\in\calP$.
For $p\in[0,1]^\N$
and $b\in\set{-1,1}^\N$,
we say that 
$$p'=p'(p,b)\in[0,1]^\N$$
is a {\em $b$-reflection of $p$}
about $\frac12$
if 
$$p'_j=b_j\paren{\p[j]-\frac12}+\frac12
,
\qquad
j\in\N.
$$
We say that
$\calP\subset[0,1]^\N$ is {\em 
strongly
symmetric about $\frac12$} if 
$
p'(p,b)\in\calP$
for all $p\in\calP$ and $b\in\set{-1,1}^\N$.

The family $\lgc\subset[0,\frac12]^\N$
was defined in \citet{CohenK23} as
the one learnable by the EME $\hat p$,
and characterized therein as consisting
precisely of those $p\in[0,\frac12]^\N$
for which $T(p)<\infty$. Since in this work
we do not restrict the range of $p$ to $[0,\frac12]$,
we define 
\beqn
\dot{\lgc}&:=&\set{p\in[0,1]^\N: 
T(p)<\infty
},
\eeqn
for $T$ as defined in \eqref{eq:Tnew}.
It is straightforward to extend the arguments
of \citet{CohenK23}
to show that
$\dot{\lgc}$
consists precisely of those $p\in[0,1]^\N$
for which the EME $\hat p$
yields $\Delta_n\to0$.

\section{Main Results}

Our first result may be informally summarized
thus: ``morally'' speaking, 
$\lgc$
is the largest 
family
that is learnable by any
fixed estimator.

\begin{theorem}[expanding $\lgc$]
\label{thm:EMELGC}

Suppose that $\calP\subset[0,1]^\N$
defines a family of product distributions
as in \eqref{eq:mup} and furthermore
\begin{enumerate}
    \item 
    $\calP$ is decaying
    \item $\calP$ is strongly symmetric about $\frac12$
    \item $\calP$ is learnable.
\end{enumerate}
Then $\calP\subseteq\dot{\lgc}$.

\hide{old:

Let $\lgc$ be the set of product distributions for which the Maximum Likelihood Estimator (EME) $\hat{p}_n$ satisfies $\Delta_n(p) \to 0$. There does not exist an estimator $\tilde p$ such that $\widetilde{\lgc}$, the set for which $\tilde{\Delta}_n(p) \to 0$, is strictly larger than $\lgc$ while satisfying the following conditions:

Formally, for any estimator $\tilde p$, defined as a mapping from $(\{0,1\}^\N)^n$ to $[0,1]^\N$, the set
\[
\widetilde{\lgc} := \left\{ p \in [0,1]^\N : \tilde{\Delta}_n(p) \to 0 \right\}
\]
cannot be strictly larger than $\lgc$, meaning that:
\[
\lgc \subseteq \widetilde{\lgc} \quad \text{and} \quad \widetilde{\lgc} \setminus \lgc \neq \emptyset,
\]
if it satisfies the following conditions:
\begin{enumerate}
    \item Symmetry: For any sign vector $b \in \{-1, 1\}^\N$, $p \in \widetilde{\lgc} \iff \frac{1}{2} - \frac{b}{2} - b p \in \widetilde{\lgc}$.
    \item Convergence of $\tilde{\Delta}_n(p)$: For all $p \in \widetilde{\lgc}$, the estimator $\tilde p$ satisfies $\tilde{\Delta}_n(p) \to 0$ as $n \to \infty$.
    \item Decay of $p$: For all $p \in \widetilde{\lgc}$, $\lim\sup_{j\in \N} p(j) \wedge (1 - p(j)) = 0$. 
\end{enumerate}

Here, $\tilde{\Delta}_n(p)$ is defined as:
\[
\tilde{\Delta}_n(p) := \E \sup_{j \in \N} |\tilde p(j) - p(j)|.
\]
}
\end{theorem}
\paragraph{Remark.}
Strong symmetry about $\frac12$
forces the sequences in $\calP$ to be ``generic''
and prevents the learner from beating the EME
by exploiting some special structure. Note that this condition is very much absent in Theorem~\ref{thm:relax}, where indeed such exploits become possible.

\hide{
\paragraph{Remark.}
The exceptional status of the sequence
$q=    (\frac12,\frac12,\ldots)$
is due to its unique property of being invariant
under $b$-reflections about $\frac12$.
As Theorem~\ref{thm:relax} shows,
$\lgc\cup\set{q}$ is learnable even though
$q\notin\dot{\lgc}$.
}




\newcommand{\ninfty}[1]{\lVert #1\rVert_{\infty}}
\newcommand{\DKL}{D_{\textup{KL}}}
\newcommand{\DKLb}{h}
\newcommand{\mmc}{\calP}


Having established that (modulo pathologies) $\lgc$ is the largest learnable family, we next show that
the EME 
is 
nearly
minimax-optimal for this family.

\begin{theorem}[Minimax bound]
\label{thm:minimax}
There exist universal constants $c,c',C>0$ such that the following holds.
For 
$n\in\N$
and $s,t>0$
satisfying
$\frac{c'\log n}{n}\leq \frac{s}{t} \leq 
\mathe^{-1}
$, 
let
\begin{align*}
\mmc_{s,t}:=
\left\{
p\in[0,1]^\N: S(p)\leq s \,\wedge\, T(p)\leq t \right\}.
\end{align*}
Then,
whenever
$\frac{c'\log n}{n}\leq\frac{s}{t}\leq 
\mathe^{-1}
$ and 
\begin{align*}
    n\geq 
    \frac{
    \frac{t^2}{Cs}\log\frac{t}{s}
    }
    {
    \frac{t}{s}\log\frac{t}{s}\cdot\mathe^{-\frac{t\log\frac{t}{s}}{\log 2}}-1},
\end{align*}
we have
\begin{align}
\label{eq:minmax_bound}
\inf_{\tpn} \sup_{p \in \mmc_{s,t}} \E \sup_{j \in \N} |\tpn[j] - \p[j]| \geq  
1 \wedge \paren{c
\sqrt{\frac{s}{n}}\,\vee\,Q(t,s)\cdot\frac{t}{n}}
,
\end{align}
where the infimum is over all estimators $\tpn$ 
that are based on $n$  i.i.d.\ samples 
drawn from $p$,
and
\begin{align*}
Q(t,s)=C\left(1 + \frac{\log\frac{t}{s}}{\log\log\frac{t}{s}}\right)^{-1}.
\end{align*}
\end{theorem}
\paragraph{Remark.}
The 
logarithmic factor and restrictions on the range
of $n$ are likely artifacts of the argument,
which we kept streamlined for space and readability. We look forward to removing both
in the extended version.


Finally, we show that if the learner is allowed
to ``cheat'' by exploiting the information
contained in the infinitely many bits of 
each example
$X\supr i$,
then $\lgc$ can indeed be non-trivially expanded.
Let us elaborate a bit on the nature of these exploits.
The elements of $\lgc$ have a ``generic,'' 
unstructured
flavor: knowing the 
values of $\p[j]$ for $j\in[N]$
reveals no 
useful
information regarding the remaining $j>N$;
all the learner knows is that these
must decay as some power of $j$ in order to be in $\lgc$.
On the other hand, one might consider adjoining a ``structured'' sequence to $\lgc$, such as
$p=(\frac12,\frac12,\ldots)$.
Because a single $X\supr i\sim\mu$ provides a bit drawn from
{\em each} of the $\Bernu(\p[j])$, the learner (as we show below) is able to first test whether the unknown 
sequence
has
the given structure
(in this case, 
whether it was generated by
$p\equiv\frac12$)
and if not, then reverts to the standard EME
for learning the unstructured 
sequences
in
$\lgc$.

\begin{theorem}[Relaxing decay and symmetry]
\label{thm:relax}
Define the family $
\mathrm{const}
\subset[0,1]^\N
$
by
\beq
\mathrm{const}
&:=& 
\set{(c, c, \dots):c\in[0,1]}
.
\eeq
Then 
$
\calP=
\dot{\lgc}\cup
\mathrm{const}
$ 
is learnable, meaning 
that
there exists an estimator $\tilde p$ such that $\tilde{\Delta}_n(p) \to 0$ for all $p \in 
\calP
$. 
\hide{
Furthermore, if we drop the symmetry assumption, the set
\[
\lgc \cup \{(c, c, \dots) \mid c \in [0, 1]\}.
\]
is learnable.
}
\end{theorem}
\paragraph{Remark.}
The techniques of Theorem~\ref{thm:relax}
are applicable 
considerably more broadly than just to
the family
$\calQ=
\mathrm{const}
$. 
For example,
the argument can be easily adapted
to show that 
$
\dot{\lgc}\cup\calQ
$
is learnable for any finite $\calQ$.
The following proposition shows a nontrivial $\calQ$ which is strongly symmetric and learnable, implying that the decay condition of Theorem~\ref{thm:EMELGC} is necessary.

\begin{proposition}[Relaxing decay]
\label{prp:relax}

Let $c > 0$ be a universal constant, and define the family of distributions
\begin{equation*}
\calQ = \set{ \p \in [0,1]^\N \mid \forall j \in \N, \abs{\p[j] - 1/2} \leq \frac{c}{\sqrt{j}} }.
\end{equation*}
Then, $\calQ$ is learnable.
\end{proposition}

\paragraph{Open problems.}
Two natural directions for future study are extensions of Theorems~\ref{thm:EMELGC} and \ref{thm:relax}.
For the former, it is likely that the conditions on $\calP$ are too stringent and can be significantly relaxed; in particular, requiring that $\calP$ be decaying
is quite probably unnecessary. Thus, we seek a larger 
family
$\calP'$  whose learnability implies
$\calP'\subseteq\dot{\lgc}$. 
Regarding Theorem~\ref{thm:relax},
we again
anticipate
the existence
of
considerably richer families $\calQ$ for which $\dot{\lgc}\cup\calQ$ is learnable.
One such family is proposed in the conjecture below.

\paragraph{Conjecture.}
Let $\calQ\subset[0,1]^\N$ be a 
countable
family of sequences with
the following property:
for each $q,q'\in\calQ$,
there is an $\eps>0$ 
and an infinite $J\subset\N$
such that $|q_j-q'_j|>\eps$
for all $j\in J$.
Then $\lgc\cup\calQ$ is learnable.


\section{Proofs}


\subsection{Proof of Theorem \ref{thm:EMELGC}}
\hide{
\begin{lemma}[$\lgc$ is permutation-invariant]
For 
all permutations $\sigma:\N\to\N$,
we have
$p\in\lgc
\implies
\sigma(p)\in\lgc$,
where
$q=\sigma(p)$
is the sequence given by $q_j=p_{\sigma(j)}$.
\end{lemma}
}
Assume, for the sake of contradiction, that there exists an estimator $\tilde{p}$ and a family $\calP\subset [0,1]^\N$ satisfying the conditions of the theorem, such that $\calP$ is learnable by $\tilde{p}$ but there exists 
a
$p^*\in\calP
\setminus
\dot{\lgc}$.
Based on $p^*$,
we will construct a family $\calP^*\subset\calP$
and argue that 
$\tilde{\Delta}_n(p) \to 0$
cannot hold for all $p\in\calP^*$.

Since $\calP$ is 
strongly symmetric about $\frac{1}{2}$, for any $p\in\calP$ and any sign vector $b\in\{-1,1\}^\N$, the $b$-reflection $p'(p,b)$ defined by $p'_j = b_j(\p[j] - \frac{1}{2}) + \frac{1}{2}$ also belongs to $\calP$.

Consider the following randomized experiment:
\begin{itemize}
    \item Let $Y = (Y_j)_{j\in\N}$ be a sequence of independent Rademacher random variables, i.e., $\P(Y_j = 1) = \P(Y_j = -1) = \frac{1}{2}$.
    \item Define $p^{(Y)} \in [0,1]^\N$ as the $Y$-reflection of $p^*$ about $\frac12$, i.e. $\p[j]^{(Y)} = Y_j(\p[j]^* - \frac{1}{2}) + \frac{1}{2}$.
    \item Generate $n$ independent draws $X^{(1)}, \dots, X^{(n)} \in \{0,1\}^\N$ from the product distribution $\mu(p^{(Y)})$
    as in \eqref{eq:mup}.
\end{itemize}

The assumption that
\(\minp[j]^*
\in[0,\frac14]^\N
\)
incurs no loss of generality,
since the decay condition implies that
$\minp[j]^*\le\frac14$ will hold for all sufficiently large $j$. We can ignore \(\minp[j] \in \left(\tfrac{1}{4}, \tfrac{1}{2}\right]\) since this would only decrease the estimation error \(\tilde{\Delta}_n(p)\).


We follow the standard reduction 
from the harder problem of estimating $p\supr Y$
to the easier problem of recovering 
the sign vector \(y \in \{-1,1\}^\N\) that defines the \(y\)-reflection \(p^{(Y)}\). 
By the Neyman-Pearson lemma,
an optimal
estimator \(\hat{y}\) 
is one
that minimizes the posterior probability of error, i.e.,
\[
\hat{y} = \argmin_{y \in \{-1,1\}^\N} \P(Y \neq y \mid X = x),
\]
where \(Y\) is the random sign vector and \(X = \paren{X^{(1)}, \dots, X^{(n)}}\) denotes the observed data.

Now
\begin{align*}
1-\P&(Y \neq y \mid X = x)
\\ =&
\P(Y = y \mid X = x)
\\
=&
\P(Y_1=y_1\mid X = x) \cdot
\\&
\P(Y_2=y_2\mid X = x,Y_1=y_1)\cdot
\\
&
\P(Y_3=y_3\mid X = x,Y_1=y_1,Y_2=y_2)\cdot\ldots
.
\end{align*}
Since the $Y_j$
are mutually independent,
each of the factors
above has the simpler form
\begin{align*}
\P(Y_k=y_k\mid 
X = x,Y_1=y_1,\ldots,Y_{k-1}=y_{k-1}
)
\\=
\P(Y_k=y_k\mid 
X=x
).
\end{align*}

%
%
We conclude that
the events
$
E_j=\set{Y_j\neq y_j \mid X}
$
are
mutually independent.
Thus
\begin{align*} 
\P&(Y \neq y \mid X = x)
\\ &=
\P\paren{\bigcup_{j\in\N}
E_j
} 
\\&=
\lim_{N\to\infty}
\P\paren{\bigcup_{j=1}^N
E_j
} 
\\&=
\lim_{N\to\infty}
\alpha_N\paren{\P(E_1),\P(E_2),\ldots,\P(E_N)},
\end{align*}

where the second equality holds
by 
regularity 
of probability measures
\citep[Theorem 17.10]{kechris:95},
and $\alpha_N:[0,1]^N\to[0,1]$
is the {\em inclusion-exclusion}
function defined inductively by $\alpha_1(x)=x$
and 
\begin{align*}
\alpha_{N+1}&(x_1,x_2,\ldots,x_{N},x_{N+1})
 \\ &=
x_{N+1}+(1-x_{N+1})
\alpha_{N}(x_1,x_2,\ldots,x_{N}).
\end{align*}
By \citet[Lemma 4.2]{kontorovich12},
$\alpha_N$ is monotonically increasing in each argument.
Hence, the optimal estimator may minimize each $\P(E_j)$
individually --- 
\hide{
Expanding this yields
\begin{align*}
\hat{y} &= \argmin_{y \in \{-1,1\}^\N} \P(Y \neq y \mid X = x) \\
&= \argmin_{y \in \{-1,1\}^\N} \prod_{j \in \N} \P(Y_j \neq y_j \mid X_j = x_j) \\
&= \left( \argmin_{y_j \in \{-1,1\}} \P(Y_j \neq y_j \mid X_j = x_j) \right)_{j \in \N}.
\end{align*}
This shows that the optimal estimator for each \(Y_j\) is the one that minimizes the individual posterior error probabilities. 
}and 
so
we may define \( A_j \) as the estimator for the \( j \)-th coordinate, where \( A_j: \{0,1\}^{\N \times n} \to [0,1] \) is any mapping from the \( j \)-th row of the data matrix to an estimate of \( \p[j]^{(Y)} \).
Let \( B_j \) be the event that \( A_j \) and \( \p[j]^{(Y)} \) belong to different intervals, i.e.,
either
$A_j\in[0,\frac12)$ and 
$\p[j]^{(Y)}\in(\frac12,1]$
or vice versa.
To establish a lower bound on the error of any estimator, consider the minimax risk:
\begin{align*}
 &\inf_A  \sup_{p \in \calP} \E_X \sup_{j \in \N} |A_j - \p[j]| 
\\&\geq \inf_A \E_{Y} \E_X \sup_{j \in \N} |A_j - \p[j]^{(Y)}| \\
&\geq \inf_A \E_{Y} \E_X \sup_{j \in \N} \mathbf{1}\{B_j\} \left|\frac{1}{2} - \p[j]^{(Y)} \right| \\
&\geq 
\frac{1}{4}
\inf_A \E_{Y} \E_X \sup_{j \in \N} \mathbf{1}\{B_j\}  \\
&= \frac{1}{4} \inf_A \P_{Y,X}\left(\bigcup_{j \in \N } B_j\right) \\
&= \frac{1}{4} \inf_A \int_{x \in \{0,1\}^{\N \times n}} \P\left(\bigcup_{j \in \N } B_j \mid X = x \right) dP_X(x) \\
& \geq \frac{1}{4}  \int_{x \in \{0,1\}^{\N \times n}} \min_{\hat{y} \in \{-1, 1\}^\N} \P(Y \neq \hat{y} \mid X = x) dP_X(x).
\end{align*}
By the Neyman-Pearson lemma,
the optimal choice of $\hat{y}_j$ is according to the majority vote\footnote{
The issue of optimally breaking ties or allowing randomized decision rules
is
somewhat delicate and is exhaustively addressed
in \citet[Eq. (2.7)]{KonPin2019}.
In our setting, these do not affect the probability of error.
}
of the vector $(X^{(1)}_j, \dots, X^{(n)}_j)$. In the event that $(X^{(1)}_j, \dots, X^{(n)}_j)=(1,1,\dots,1)$, but $1-\minp[j]^* \neq \p[j]^{(Y)}$, the estimator makes a mistake. 
The probability of such an event, conditioned on the other random variables $X_{j'}, Y_{j'}$ where $j'\neq j$, is exactly $\frac{1}{2}(\minp[j]^*)^n$.
Since we assumed $\p^* \notin \lgc$ and thus $\minp^* \notin \lgc$, we have $T(\minp^*) = \infty$. 
Since $\minp[j]^* \to0$,
we may assume without loss of generality that
it is decreasing monotonically.
By \citet[Lemma 3]{CohenK23}, 
it follows
that for all $n \in \N$,
we have
\[
\sum_{j=1}^\infty (\minp[j]^*)^n = \infty.
\]
Since the events of 
$\hat{y}_j$ being wrong are 
mutually
independent, 
the second Borel–Cantelli lemma
implies that
almost surely
at least one of them will occur.
It follows that $
\liminf_{n\to\infty}
\tilde\Delta_n (
p\supr Y
) \geq \frac{1}{4}$, contradicting the learnability assumption.
\QED




\subsection{Proof of Theorem \ref{thm:minimax}}
We reduce the minimax lower bound problem to one over a finite set of hypotheses.
For $2\leq J\in\N$ and $0\leq q \leq q' \leq 1/2$ to be chosen below, we consider $J+1$ profiles $p\supr{k}\in[0,\frac 1 2]^\N$
 for $k\in [J+1]$.
For $k=J+1$ we take the step profile
\begin{align*}
    \p[j]\supr{J+1}=
    \begin{cases}
        q, &  j\in [J+1],
        \\
        0, & j>J+1,
    \end{cases}
\end{align*}
and for $1\leq k\leq J$ we take the same step profile but with an additional bump at position $k$,
\begin{align*}
    \p[j]\supr{k}=
    \begin{cases}
        q, & j\in [J+1] \text{ and } j\neq k,
        \\
        q', & j=k,
        \\
        0, & j>J+1.
    \end{cases}
\end{align*}

The construction of these profiles is illustrated in Figure~\ref{fig:profiles}.

\begin{figure}[ht!]
    \centering
    \begin{tikzpicture}[scale=0.8]
        \draw[->] (0,0) -- (9,0) node[below] {$j$};
        \draw[->] (0,0) -- (0,3.5) node[above] {\(\p[j]\supr{k}\)};

        \draw (1,0) node[below] {$1$};
        \draw (1.75,0) node[below] {$2$};
        \draw (2.5,0) node[below] {$\dots$};
        \draw (3.25,0) node[below] {$k$};
        \draw (4,0) node[below] {$\dots$};
        \draw (4.75,0) node[below] {$J$};
        \draw (5.5,0) node[below] {\rotatebox{90}{$J+1$}};
        \draw (6.25,0) node[below] {\rotatebox{90}{$J+2$}};
        \draw (7,0) node[below] {\rotatebox{90}{$J+3$}};

        \foreach \x/\h in {1/2, 1.75/2, 4.75/2, 5.5/2} {
            \draw[fill=blue!50] (\x-0.15,0) rectangle (\x+0.15,\h); 
            \node[above] at (\x,\h) {$q$};
        }

        \draw[fill=red!70] (3.25-0.15,0) rectangle (3.25+0.15,3); 
        \node[above] at (3.25,3) {$q'$};

        \foreach \x in {6.25,7} {
            \draw[fill=gray!30] (\x-0.15,0) rectangle (\x+0.15,0.1); 
            \node[above] at (\x,0.1) {$0$};
        }

        \node at (7.75,-0.15) {$\dots$};
    \end{tikzpicture}

    \begin{tikzpicture}[scale=0.8]
        \draw[->] (0,0) -- (9,0) node[below] {$j$};
        \draw[->] (0,0) -- (0,3.5) node[above] {\(\p[j]\supr{J+1}\)};

        \draw (1,0) node[below] {$1$};
        \draw (1.75,0) node[below] {$2$};
        \draw (2.5,0) node[below] {$\dots$};
        \draw (4.0,0) node[below] {$\dots$};
        \draw (4.75,0) node[below] {$J$};
        \draw (5.5,0) node[below] {\rotatebox{90}{$J+1$}};
        \draw (6.25,0) node[below] {\rotatebox{90}{$J+2$}};
        \draw (7,0) node[below] {\rotatebox{90}{$J+3$}};

        \foreach \x/\h in {1/2, 1.75/2, 3.25/2, 4.75/2, 5.5/2} {
            \draw[fill=blue!50] (\x-0.15,0) rectangle (\x+0.15,\h); 
            \node[above] at (\x,\h) {$q$};
        }

        \foreach \x in {6.25,7} {
            \draw[fill=gray!30] (\x-0.15,0) rectangle (\x+0.15,0.1); 
            \node[above] at (\x,0.1) {$0$};
        }

        \node at (7.75,-0.15) {$\dots$};
    \end{tikzpicture}

    \caption{Illustration of the step profile construction for \(p\supr{k}\) (top) and the special case for \(p\supr{J+1}\) (bottom). Each bar represents the value of \(\p[j]\supr{k}\) at position \(j\). Values are shown above the bars.}
    \label{fig:profiles}
\end{figure}

Note that for all $k\neq\ell\in[J+1]$ we have
$\ninfty{p\supr{k}-p\supr\ell} = |q'-q|$.
In addition, for $k=J+1$,
\[
    S(p\supr{J+1}) = q\log(J+1)
\]
and
\[
T(p\supr{J+1}) = \frac{\log(J+1)}{\log \frac{1}{q}},
\]
and for $k\in[J]$,
\[
S(p\supr{k}) = \max\left\{q\log(J+1), q'\log 2\right\}
\]
and
\[
T(p\supr{k}) = \max\left\{\tfrac{\log(J+1)}{\log \tfrac 1 q}, \tfrac{\log2}{\log \tfrac{1}{ q'}}\right\}.
\]
Given $s\leq 
\frac{t}{\mathe}
$ as in the Theorem statement, 
we 
choose
$q\in[0,\frac12]$
and 
$J$
as
\begin{align*}
    q=\frac{1}{\frac{t}{s}\log \frac{t}{s}}
    \qquad \text{and} \qquad
    \log(J+1) = t\log\frac{t}{s}.
\end{align*}
Then
\begin{align}
\label{eq:S_s}
S(p\supr{J+1})
&= q\log(J+1) = s,
\\[3pt]
\label{eq:T_t}
T(p\supr{J+1})
&= \frac{\log(J+1)}{\log\frac 1 q} = t 
\cdot
\left(1 + \frac{\log\frac{t}{s}}{\log\log\frac{t}{s}}\right)^{-1}
 \leq t.
\end{align}
Below we set $q'\leq 1/2$ such that $q'\geq q$ and
for all $k\in[J]$,
\begin{align}
\label{eq:in_class_cond}
S(p\supr{k})\leq S(p\supr{J+1})
= s
\nonumber
\\
\quad\text{and}\quad 
T(p\supr{k})&\leq T(p\supr{J+1}) \leq t.
\end{align}
Thus, $p\supr{k}\in \mmc_{s,t}$ for all $k\in[J+1]$ and 
\begin{align}
\label{eq:red_to_fin}
\inf_{\tpn} & \sup_{p \in \mmc_{s,t}} \E \sup_{j \in \N} |\tpn[j] - \p[j]| 
\nonumber \\
&\geq
\inf_{\tpn}
\max_{k \in [J+1]}
\E_{X_n \sim \mu\supr{k,n}}
    \ninfty{\tpn(X_n)-p\supr{k}}
.
\end{align}
To lower bound the right-hand side of \eqref{eq:red_to_fin} we apply the generalized Fano method. 
For $k\in[J+1]$, let $\mu\supr{k}
=\mu(p\supr k)
$ 
be the product measure over $\set{0,1}^\N$
as defined in \eqref{eq:mup}
and
note that $\E_{X\sim\mu\supr{k}}\{X\}=p
\supr{k}$.
We denote by 
$\mu\supr{k,n}$ 
the product measure of $n$ independent copies of $X\sim\mu\supr{k}$.
We invoke Lemma \ref{lem:fano} with the $J+1$ measures $(\nu_1,\dots,\nu_{J+1})=
(
\mu\supr{1,n},\ldots,\mu\supr{J+1,n}
)$, the distance function $\rho=\ninfty{\cdot}$, and the parameters $\theta(\mu\supr{k,n}) = \E_{X\sim\mu\supr{k}}\{X\} = p\supr{k}$ for $k\in[J+1]$.
Note that $    \rho(\theta(\mu\supr{k,n}),\theta(\mu_
\ell^n)) = |q'-q|$ for all $k\neq\ell\in[J+1]$
and that
\begin{align*}
    \DKL(\mu\supr{k,n} \lVert \mu_\ell^n) & \leq 
    n(\DKLb(q \lVert q') + \DKLb(q' \lVert q)),
\end{align*}
where 
\begin{align*}
    \DKLb(q \lVert q') = q\log\frac{q}{q'} + (1-q)\log\frac{1-q}{1-q'}.
\end{align*}
Then
Lemma \ref{lem:fano} implies
\begin{align}
\nonumber
\inf_{\tpn} &
\max_{k \in [J+1]}
\E_{X_n \sim \mu\supr{k,n}}
    \ninfty{\tpn(X_n)-p\supr{k}}
     \\ & \geq
    \frac{q'-q}{2} \paren{1 - \paren{\frac{n(\DKLb(q \Vert q') + \DKLb(q' \Vert q)) + \log 2}{\log(J+1)}}}.
    \label{eq:lb_2hpq}
\end{align}
We now fix $q'(q)=q'(q,n,J)\geq q$ to be the solution to the equation
\begin{align}
\label{eq:hqp_equal}
\DKLb(q \Vert q'(q)) +  \DKLb(q'(q) \Vert q) = \frac{\log(J+1)}{2Cn}.
\end{align}
Below we verify that 
\eqref{eq:in_class_cond}
indeed holds with this choice of $q'(q)$.
Substituting \eqref{eq:hqp_equal} into \eqref{eq:lb_2hpq}, we obtain the lower bound
\begin{align*}
\frac{q'(q)-q}{2} \paren{1 - \paren{\frac{\log(J+1)/C + \log 2)}{\log(J+1)}}}  
\\
\geq 
\frac{q'(q)-q}{2}\paren{1-\frac{1}{C}- \frac{\log 2}{\log 3}}
\\
\geq
\frac{q'(q)-q}{8},
\end{align*}
for an appropriate value of the constant $C>0$.

We analyze 
$q'(q)-q$
for $q'(q)$ satisfying \eqref{eq:hqp_equal} as in \citet{BlanchardCK24} and consider two regimes for $\DKLb(q\Vert q')+\DKLb(q'\Vert q)$.
For any $0< q\leq q' \leq \frac{1}{2}$, 
we have 
\begin{align}
\label{eq:h_ineq_s}
\frac{(q'-q)^2}{q'}    \leq \DKLb(q \Vert q') +\DKLb(q' \Vert q)  \leq \frac{2(q'-q)^2}{q}.
\end{align}
So, by the right inequality in \eqref{eq:h_ineq_s},
\begin{align}
\label{eq:qp_S_sg}
     q'(q)-q & \geq 
     \sqrt{\frac{q\left(\DKLb(q\Vert q'(q))+\DKLb(q'(q)\Vert q)\right)}{2}}
     \nonumber \\ &=
     \sqrt{\frac{q\log(J+1)}{4Cn}} 
     \nonumber \\ &= \sqrt{\frac{S(p_{J+1})}{4Cn}} 
     \nonumber \\ &= \sqrt{\frac{s}{4Cn}}.
\end{align}
In addition, by the left inequality in \eqref{eq:h_ineq_s},
\begin{align*}
    q' &\leq q + \sqrt{\frac{q'\log(J+1)}{2Cn}}
    \\ &\leq \sqrt{q'}
    \left(\sqrt{q} + \sqrt{\frac{\log(J+1)}{2Cn}}\right),
\end{align*}
which implies
\begin{align*}
    q' & \leq \left(\sqrt{q} + \sqrt{\frac{\log(J+1)}{2Cn}}\right)^2 
    \\&= q + 2\sqrt{\frac{q\log(J+1)}{2Cn}} + \frac{\log(J+1)}{2Cn}
    \\ &=
    q\left(1+ 2\sqrt{\frac{\log(J+1)}{2Cnq}} + \frac{\log(J+1)}{2Cnq}\right).
\end{align*}
Since by assumption $Cqn = \frac{Cn}{\frac{t}{s}\log(\frac{t}{s})} \geq c'\log 2$,
we have that for a sufficiently large constant $c'$,
\begin{align*}
    q' \leq \frac{q\log(J+1)}{\log 2}.
\end{align*}
This verifies \eqref{eq:in_class_cond}
and establishes the term $c\sqrt{\frac{s}{n}}$ in \eqref{eq:minmax_bound}.


Next, we assume $\frac{t}{n}\geq  c\sqrt{\frac{ s}{n}}$.
For any $0\leq q\leq q' \leq 1/2$ we have
$
\DKLb(q \Vert q') \leq \DKLb(q' \Vert q)$,
and for $q'\geq 9q$, it holds that (see, e.g., 
\citet{BlanchardCK24})
\begin{align}
\label{eq:h_ineq_t}
\frac{1}{2} \leq
    \frac{\DKLb(q \lVert q') + \DKLb(q' \lVert q) }{q\cdot\tfrac{q'-q}{q}\log\tfrac{q'-q}{q}}
    \leq
    4.
\end{align}
For $z\geq \mathe$, the solution $x$ to the equation $x\log x=z$ satisfies $x\geq \frac{z}{\log z}$ 
\citep{lambertW}.
Taking $c>0$ sufficiently large such that
\begin{align*}
z = \frac{\log(J+1)}{2Cqn} = 
\frac{t^2}{2Cns}\log^2\left(\tfrac{t}{s}\right) \geq 
\frac{c^2}{2C}
\geq \mathe,
\end{align*}
we have that $q'(q)$ satisfies \eqref{eq:hqp_equal} if
$q'(q)-q\geq 8q$ and
\begin{align*}
    \frac{q'(q)-q}{q}\geq 
    \frac{\frac{\log(J+1)}{2Cqn}}{\log \frac{\log(J+1)}{2Cqn}};
\end{align*}
namely,
\begin{align*}
\label{eq:qp_T_se}
q'(q)-q&\geq 
    \frac{\log(J+1)}{2Cn\log \paren{\frac{\log(J+1)}{2Cqn}}}
    \\&= \frac{T(p_{J+1})}{2Cn} 
    \cdot
    \frac{\log \frac{1}{q}}{\log \frac{\log(J+1)}{2Cqn}}
    \\&= \frac{T(p_{J+1})}{2Cn} 
    \cdot
    \frac{\log \frac{1}{q}}{\log \frac{s}{2Cn q^2}}
    \\&\geq 
    \frac{T(p_{J+1})}{4Cn}
    \\&=
    \frac{t}{4Cn}
     \cdot\left(1 + \frac{\log\frac{t}{s}}{\log\log\frac{t}{s}}\right)^{-1}
    ,
\end{align*}
where in the last inequality we used the fact that $\frac{s}{2Cn}\leq 1$.
Lastly, we verify that $q'(q)$ is such that \eqref{eq:in_class_cond} holds, namely,
$   \frac{\log\frac{1}{q}}{\log\frac{1}{q'(q)}}\leq \frac{\log(J+1)}{\log 2}.$
The left inequality in \eqref{eq:h_ineq_t} implies that $q'(q)\leq q + C'\frac{t}{n}$ for some constant $C'$. Putting this and $q=\frac{1}{\frac{t}{s}\log\frac{t}{s}}$ and $\log(J+1)=t\log\frac{t}{s}$, we have that \eqref{eq:in_class_cond} holds if 
\begin{align*}
\frac{\log(\frac{t}{s}\log\frac{t}{s})}{\log(\frac{t}{s}\log\frac{t}{s})-\log(1 +  \frac{t^2}{Cns}\log\frac{t}{s})}\leq \frac{t\log\frac{t}{s}}{\log 2}.
\end{align*}
This is satisfied when
\begin{align*}
    n\geq 
    \frac{
    \frac{t^2}{Cs}\log\frac{t}{s}
    }
    {
    \frac{t}{s}\log\frac{t}{s}\cdot\mathe^{-\frac{t\log\frac{t}{s}}{\log 2}}-1}.
\end{align*}
Finally, we consider the case where $t\geq n$. We repeat the arguments in the proof of Theorem~\ref{thm:EMELGC} to show that in this case the minimax rate is bounded from below by a constant. Taking any $p^* \in \mmc_{s,t}$ such that $T(p^*) \geq n$, and assuming without loss of generality that $p^*$ is non-increasing, 
let $j'$ be such that
\begin{align*}
    T(p^*) \geq \frac{\log(1+j')}{\log (1/\minp[j']^*)}.
\end{align*}
Then
\begin{align*}
\sum_{j=1}^\infty (\minp[j]^*)^n 
& \geq \sum_{j=1}^\infty (\minp[j]^*)^{T(p)}
\\& \geq
\sum_{j=1}^{j'}(\minp[j]^*)^{\frac{\log(1+j')}{\log (1/(\minp[j]^*))}}
\\& \geq 
j'(\minp[j']^*)^{\frac{\log(1+j')}{\log (1/(\minp[j']^*))}}
\\& =
 \frac{j'}{1+j'}
\geq \frac{1}{2}.
\end{align*}
As in the proof of Theorem \ref{thm:EMELGC}, applying
Lemma~\ref{lem:ramon}
with
\begin{align*}
&A_j = \\ &\set{(X^{(1)}_j, \dots, X^{(n)}_j)=(1,1,\dots,1),
\text{ but }
1-\minp[j]^* \neq \p[j]^{(Y)}},
\end{align*}

where $Y_j\sim \Bernu(1/2)$ and $\p[j]\supr{Y_j}=Y_j\minp[j]^*+(1-Y_j)(1-\minp[j]^*)$, we get that the minimax rate is lower bounded by a universal constant.
\QED

\subsection{Proof of Theorem \ref{thm:relax}
}

We aim to prove that the family
\[
\calP \eqdef \dot{\lgc} \cup \set{(c, c, \dots) : c \in [0,1]}
\]
is learnable by an estimator $\tilde p_n$. 
Choose some $p \in \calP$. The general strategy is to construct an estimator that can distinguish between cases where $T(p) = \infty$ and cases where $T(p) < \infty$, 
based on
the sample.

\paragraph{Step 1: Testing if $T(p) = \infty$.}
We begin by defining a test $\Phi$ to check whether $T(p) = \infty$. The idea is to check the first half of the sequence $X\supr1_j,X\supr2_j,\ldots,X\supr n_j$ are all ones and the second half are all zeros. Formally, the test function is defined as:
\[
\Phi(X\supr1,X\supr2,\ldots,X\supr n) = \mathbf{1}\paren{ \limsup_{j \to \infty} E_j},
\]
where we define the event
\[
E_j \eqdef  
\set{X\supr i_j =0 \text{ for } i \leq \frac{n}{2} \text{ and } X\supr i_j =1 \text{ for } i > \frac{n}{2}}.
\]
Note that, for each $j$, we have
\[
\P \paren{E_j} = \p[j]^{\floor{n/2}} \paren{1-\p[j]}^{\ceil{n/2}},
\]
and because $\set{E_j}_j$ are independent, by the two Borell-Cantelli lemmas, $\Phi = 1$ almost surely if and only if
\[
\sum_{j=1}^\infty \p[j]^{\floor{n/2}} \paren{1-\p[j]}^{\ceil{n/2}} = \infty.
\]
The above sum can be estimated by the following sums,
\beqn
\label{eq:series_bounds}
\sum_{j=1}^\infty \minp[j]^{n} 
\leq
\sum_{j=1}^\infty \p[j]^{\floor{n/2}} \paren{1-\p[j]}^{\ceil{n/2}}
\leq
\sum_{j=1}^\infty \minp[j]^{\floor{n/2}} .
\eeqn

\paragraph{Step 2: Consistency of the Test.}
We now show that the test $\Phi$ is consistent. First, assume $T(p) = \infty$, which means $T(\dot p_{\downarrow0}) = \infty$, then by 
 \citet[Lemma~3]{CohenK23}, we have $\sum_{j=1}^\infty \minp[j]^n = \infty$ for all $n$, then by \eqref{eq:series_bounds} we have $\Phi=1$ almost surely.

On the other hand, if $T(p) < \infty$, again by \citet[Lemma~3]{CohenK23}, we have $\sum_{j=1}^\infty \minp[j]^n < \infty$ for large enough $n$, which means that for large enough $n$ we have $\Phi = 0$ almost surely, as before.

\paragraph{Step 3: Defining the Estimator.}
Once the test $\Phi$ has been applied, we define the estimator $\tilde{p}_n$ as follows:
\[
\tilde{p}_n(j) =
\begin{cases}
\frac{1}{n} \sum_{i=1}^n \hat{p}_n(i), & \text{if } \Phi(X_1, \dots, X_n) = 1, \\
\hat{p}_n(j), & \text{otherwise}.
\end{cases}
\]
In words, if the test $\Phi$ indicates that $T(p) = \infty$, we use the average of all $\hat{p}_n(j)$ (as this is consistent with the assumption that $p$ is a constant sequence). If $\Phi$ indicates that $T(p) < \infty$, we use the EME $\hat{p}_n(j)$ directly.

\paragraph{Step 4: Consistency of the Estimator.}
We now verify that the estimator $\tilde{p}_n$ is consistent.

If $p \in \lgc$, then $T(p) < \infty$, and the test $\Phi$ will eventually return $0$. In this case, the estimator $\tilde{p}_n(j)$ is simply the EME, which is known to be consistent for all $p \in \lgc$. Therefore, $\tilde{p}_n \to p$ in $\ell_\infty$ as $n \to \infty$.

If $p = (c, c, \dots)$ for some constant $c \in [0, 1]$, then $T(p) = \infty$, and the test $\Phi$ will always return $1$. In this case, the estimator $\tilde{p}_n(j)$ is the average of all $\hat{p}_n(j)$, which, by the law of large numbers, will converge to $c$. Thus, $\tilde{p}_n(j) \to c$ as $n \to \infty$.

\paragraph{Conclusion.}
The estimator $\tilde{p}_n$ correctly learns all distributions in the family $\dot{\lgc} \cup \set{(c, c, \dots) : c \in [0,1]}$, completing the proof.

\QED

\subsection{Proof of Proposition~\ref{prp:relax}}

Define the estimator $\tilde{p}_n(j)$ as follows:
\begin{itemize}
    \item For indices $j \leq k(n)$, where $k(n)$ is chosen later, estimate $\p[j]$ using the empirical mean:
    \begin{align*}
        \tilde{p}_n(j) = \frac{1}{n} \sum_{i=1}^{n} X_j^{(i)},
    \end{align*}
    where $X_j^{(i)} \sim \Bernu(\p[j])$.
    \item For indices $j > k(n)$, set $\tilde{p}_n(j) = 1/2$.
\end{itemize}

Using Chernoff's bound, for any $\epsilon > 0$,
\begin{align*}
    \PR{\abs{\tilde{p}_n(j) - \p[j]} \geq \epsilon} \leq 2 \exp(-2n\epsilon^2).
\end{align*}
Applying the union bound over $j \leq k(n)$,
\begin{align*}
    \PR{ \sup_{j \leq k(n)} \abs{\tilde{p}_n(j) - \p[j]} \geq \epsilon } \leq 2 k(n) \exp(-2 n \epsilon^2).
\end{align*}
Choosing
\begin{align*}
    \epsilon_n = O\left(\sqrt{\frac{\log k(n)}{n}}\right),
\end{align*}
ensures this probability vanishes. Thus, with high probability,
\begin{align*}
    \sup_{j \leq k(n)} \abs{\tilde{p}_n(j) - \p[j]} = O\left(\sqrt{\frac{\log k(n)}{n}}\right).
\end{align*}

For $j > k(n)$, since $\p[j]$ is approximated by $1/2$:
\begin{align*}
    \sup_{j > k(n)} \abs{\p[j] - 1/2} \leq \sup_{j > k(n)} \frac{c}{\sqrt{j}} = O\left(\frac{1}{\sqrt{k(n)}}\right).
\end{align*}
Choosing $k(n) = \Theta(n)$ ensures this term vanishes.

Thus, combining both terms,
\begin{align*}
    \E \norm{\tilde{p}_n(j) - \p}_\infty = O\left(\sqrt{\frac{\log k(n)}{n}}\right) + O\left(\frac{1}{\sqrt{k(n)}}\right),
\end{align*}
which converges to zero as $n \to \infty$. Therefore, $\calQ$ is learnable.

\QED

\subsection{Auxiliary lemmas}

\begin{lemma}[\citet{yu1997assouad}]
\label{lem:fano}
    For 
    \(r \geq 2\), let 
        \(\nu_1, \nu_2, ... ,\nu_r\)
    be a collection of $r$
    probability measures
    with some parameter of interest
    \(\theta(\nu)\) 
    taking
    values in pseudo-metric space \( (\Theta, \rho) \)
    such that for all
    \(j \neq k \),
    \[
    \rho(\theta(\nu_j), \theta(\nu_{k}) )
    \geq
    \alpha
    \]
    and
    \[
    \DKL(\nu_j \Vert \nu_{k})
    \leq
    \beta.
    \]
    Then
    \[
    \inf_{\hat\theta}
    \max_{k \in [r]}
    \E_{Z \sim \nu_k}
    \rho(\hat\theta(Z), \theta(\nu_k) )
    \geq
    \frac{\alpha}{2} \paren{1 - \paren{\frac{\beta + \log 2}{\log r}}},
    \]
    where the infimum is over all estimators 
    \(\hat\theta:Z\mapsto\Theta\).
\end{lemma}

\begin{lemma}[\citet{van2014probability} 
Problem~5.1a]
\label{lem:ramon}
If $A_1,\ldots,A_N$
are independent events, then
\beq
(1-\mathe\inv)\sqprn{
1\wedge
\sum_{k=1}^N \P(A_k)
}
&\le&
\P\paren{
\bigcup_{k=1}^N
A_k
}.
\eeq
\end{lemma}


\bibliographystyle{plainnat}
\bibliography{refs}

\newpage
\appendix
\onecolumn
\section{Simulation Results}

To support our theoretical findings, we present two sets of simulations. The first demonstrates the tightness of the lower bound in Theorem~\ref{thm:minimax}, while the second highlights a specific setting where the simple average estimator outperforms the Empirical Mean Estimator (EME), complementing the results of Theorem~\ref{thm:relax}.

\subsection{Tightness of the Lower Bound in Theorem~\ref{thm:minimax}}
The first simulation aims to validate the theoretical bounds presented in Theorem~\ref{thm:minimax}. Specifically, we compare the empirical average supremum deviation $\Delta_n$ with the theoretical predictions for different values of $q$ (variance control parameter) and sample sizes $n$.

We consider six values of $q$: $q = 0.1$, $q = 0.2$, $q = 0.05$, $q = 0.01$, $q = 0.005$, and $q = 0.002$. For each configuration, empirical results are averaged over $J = 100$, $1000$, and $10000$ repetitions to ensure stability. The empirical deviations are plotted alongside theoretical predictions in a log-log scale to capture the decay behavior as $n$ increases.

\begin{figure}[ht!]
    \centering
    \includegraphics[width=0.9\textwidth]{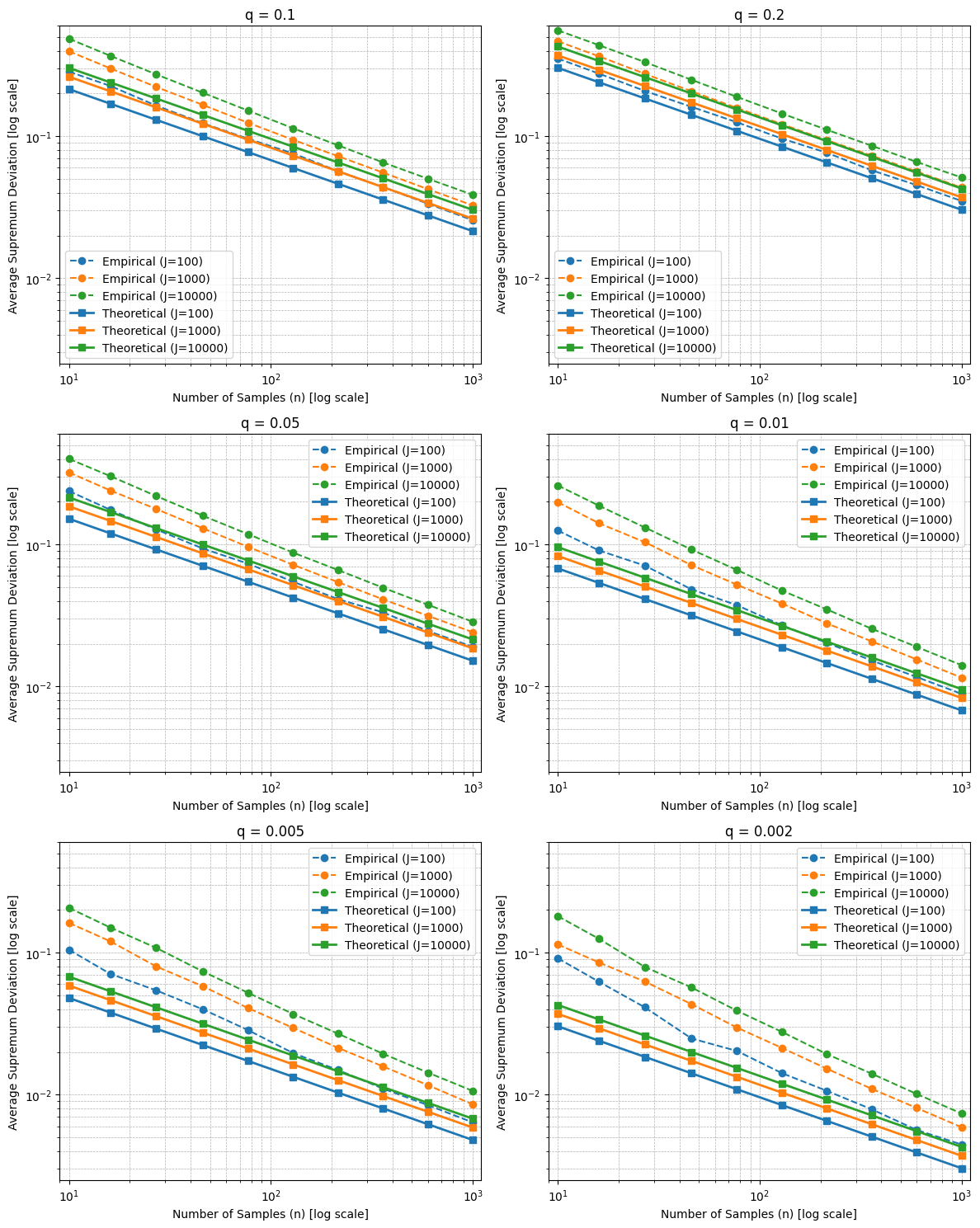}
    \caption{Average supremum deviation $\Delta_n$ as a function of sample size $n$ on a log-log scale for varying $q$ values ($q = 0.1, 0.2, 0.05, 0.01, 0.005, 0.002$). Empirical results (dashed lines) are averaged over $J = 100$, $1000$, and $10000$ repetitions and are compared to theoretical predictions (solid lines).}
    \label{fig:tightness_lb}
\end{figure}

Figure~\ref{fig:tightness_lb} shows the results. The empirical deviations (dashed lines) closely follow the theoretical bounds (solid lines), confirming the tightness of the lower bound in Theorem~\ref{thm:minimax}. As expected, larger values of $J$ lead to smoother empirical curves, emphasizing the role of averaging in reducing variance. Notably, the empirical deviations converge to the theoretical decay rate as $n$ grows.

\subsection{Performance Comparison: Theorem~\ref{thm:relax} Complement}
In the second simulation, we explore a specific setting where the simple average estimator surpasses the performance of the EME. This complements the findings of Theorem~\ref{thm:relax} by demonstrating that allowing certain structured distributions can yield better decay rates with alternative estimators.

\begin{figure}[ht!]
    \centering
    \includegraphics[width=1\textwidth]{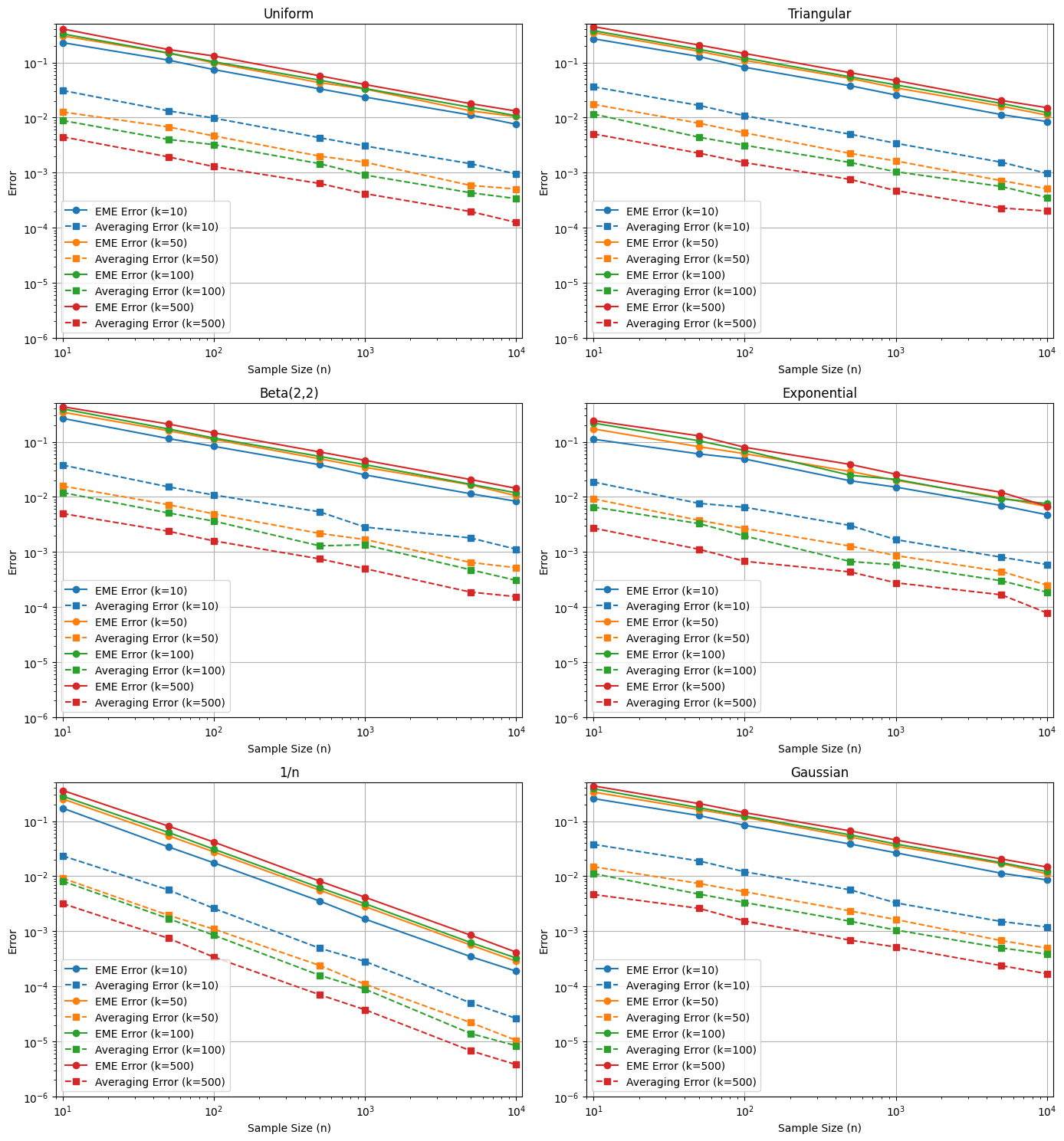}
    \caption{Error comparison between the EME and the simple average estimator for varying sample sizes $n$ under different distributions: uniform, triangular, Beta(2,2), exponential, $1/n$, and Gaussian. Results are plotted for $k \in \{10, 50, 100, 500\}$ to illustrate the effect of averaging.}
    \label{fig:avg_vs_eme}
\end{figure}

We evaluate the performance of the EME and the simple average estimator under six different distributions: uniform, triangular, Beta(2,2), exponential, $1/n$-scaled, and Gaussian. For each distribution, we vary the number of trials $k \in \{10, 50, 100, 500\}$ and compute the error as a function of the sample size $n$. The results, plotted on a log-log scale, are shown in Figure~\ref{fig:avg_vs_eme}.

The plots reveal that the simple average estimator achieves lower error rates compared to the EME across all settings as $k$ increases. This improvement is most pronounced for structured distributions like $1/n$-scaled, Beta(2,2), and Gaussian, where the averaging process effectively captures the underlying structure. These findings corroborate the theoretical insights of Theorem~\ref{thm:relax}, showcasing that the choice of estimator can significantly impact performance in specific scenarios.

\subsection{Discussion}
The results of these simulations provide strong empirical support for our theoretical findings. The first simulation confirms the tightness of the lower bound in Theorem~\ref{thm:minimax}, while the second demonstrates the practical advantages of alternative estimators, as predicted by Theorem~\ref{thm:relax}. These findings highlight the robustness and relevance of our theoretical framework for analyzing the Local Glivenko-Cantelli (LGC) class.

\end{document}